\documentclass{article}
\usepackage{a4wide}
\usepackage{graphicx}
\usepackage{amsfonts}
\usepackage{xspace}
\newcommand{\Figref}[1]{Figure \ref{#1}}

\newcommand{\figref}[1]{Figure \ref{#1}}
\newcommand{\eqnref}[1]{Equation (\ref{#1})}
\newcommand{\secref}[1]{Section \ref{#1}}
\newcommand{\ito}{It\^o }
\newcommand{\strat}{Stratonovich\xspace}

\begin{document}

\title{Effects of noise on models of spiny dendrites
%\thanks{Grants or other notes about the article that should go on the
%front page should be placed here. General acknowledgments should be
%placed at the end of the article.} 
}


\author{Emma J. Coutts         \and
        Gabriel J. Lord %etc.
}


%\address{Emma J. Coutts \at
%              \email{emmac@ma.hw.ac.uk}           %  \\
%%             \emph{Present address:} of F. Author  %  if needed
%           Gabriel J. Lord \at
%           Department of Mathematics\\
%Heriot Watt University\\
%Edinburgh\\
%EH14 4AS\\
%\email{g.j.lord@hw.ac.uk} 

\maketitle

\begin{abstract}
We study the effects of noise in two models of spiny
dendrites. Through the introduction of different types of noise to both the
Spike-diffuse-spike (SDS) and Baer-Rinzel (BR) models we investigate
the change in behaviour of the travelling wave solutions present in
the deterministic systems, as noise intensity increases. We show that
the speed of wave propagation in the SDS and BR models respectively
decreases and increases as the noise intensity in the spine heads
increases. Interestingly the discrepancy between the models does not seem to
arise from the type of active spine head dynamics employed by the
model but rather by the form of the spine density used. 
In contrast the cable is very robust to noise and as such the speed shows very
little variation from the deterministic system.  

We look at the effect of the noise interpretation used to evaluate the
stochastic integral; \ito or \strat and discuss which may be
appropriate.
%the
%approbetter option for neural models. 
%
We also show that the correlation time and length scales of
the noise can enhance propagation of travelling wave solutions where
the white noise dominates the signal and produces noise induced
phenomena. 
%\keywords{Spiny dendrite \and Stochastic wave \and \ito \strat \and space-time white \and spatial correlation length \and}
%% \PACS{PACS code1 \and PACS code2 \and more}
\end{abstract}

\section{Introduction}
\label{sec:intro}
The neuron, or nerve cell, is the building block of the mammalian
nervous system; it sends, receives and processes information that
ultimately controls functions as fundamental as our breathing and as
complex as memory. The neuron comes in many forms depending on which
area of the brain it occupies and its function but all neurons share
the same basic structure. 
%All neurons consist of the main cell body,
%or soma, an axon along which information travels to the synaptic
%terminals to be transferred onto several other connecting neurons
%through their dendritic trees. 
We are interested in models of spiny dendritic tissue and the effects of noise on these signal processing
capabilities. The dendritic tree allows a greater surface area for
synaptic connections and around $90\%$ of excitatory synapses in the
brain are made onto dendritic spines. Dendrites are typically 1-2 mm
long and the dendritic spines are small bulbous protrusions of 1-2
$\mu$m long.  Spiny dendrites occur in many regions of the brain
e.g. CA1 and CA3 pyramidal neurons in the hippocampus (important in
long term memory), basal ganglia (used in motor control and learning)
and spiny stellate neurons in the cerebral cortex (also important in
memory) \cite{Thompson}, \cite{cortex}, \cite{hippo}. The spines are
thought to be an important component in signal propagation and
computations along the dendrite and spine motility and morphology,
called spine plasticity, is thought to be an important process in
learning and memory, \cite{Yuste}, \cite{Dan}, \cite{Bonhoeffer} and
\cite{Kasai}.  In the models investigated here the plasticity can be
related to physical properties of the dendritic tissue such as the
spine stem resistance and spine density.  
The advent of the confocal and two-photon microscopy used to image the
membrane in dendrites allowed the measurement of action potentials
(APs) in dendrites and proved that action potentials can be generated
in the dendrites themselves. These techniques made it possible to
compare experimental, \cite{Hausser}, \cite{Yuste}, and theoretical
results, \cite{Vetter}, \cite{Rudolph}, that predict how voltage will
spread throughout a length of dendrite, or a branched dendritic
structure, also see for a review \cite{Segev}.  

We consider voltage spread throughout a length of spiny dendrite as a
wave propagating from the spines at the distal end, through the main
body of the dendrite to the soma without including the effect of the
soma.  This is an interesting problem as much information processing
may occur prior to the action of the soma, see review \cite{Mel}. 
We use two dendritic models that describe voltage evolution in a
length of spiny dendrite: the Baer-Rinzel (BR) model, \cite{Baer}, and
the Spike-Diffuse-Spike (SDS) model, \cite{Coombes1},
\cite{Coombes2}. Both these models couple active spines to a passive
cable through a spine density which is a constant, continuum value in
the original BR model and equally distributed point attachments in the
SDS model. We extend the BR model to include a spatially dependent
spine density which can be controlled through one parameter; the
density is equivalent to a continuum in one limit and to the point
attachment in the other. In this way we investigate the importance of
the spine stem area on the propagation of waves in the dendrite
models. 
We consider the effect of random fluctuations, or noise, in the BR and
SDS models. There are two types of noise in a neural system, intrinsic
and extrinsic noise, \cite{Gerstner}, \cite{Manwani}, (see
\cite{Faisal} for overview of noise in all levels of the central
nervous system).  Intrinsic noise is a source of noise which is always
present in the system and thermal noise is one example. Another source
of intrinsic noise, which may be considered particularly relevant in
spiny dendrites, is the process of synaptic gating since the release
of the neurotransmitter is a stochastic process. The second type of
noise, extrinsic, emanates from out with the cell itself and one
source is other, nearby neurons. 
One interpretatoin of the correlation length associated with our
spatially correlated noise is the length scale over which signals as 
inputs are transmitted. We see that there are length scales that
promote propagation.

%{\bf can we tie in somethign like : We also consider noise as
%  signal from neurons via the synapse. Conceptually spatial
%  correlation relates to connectons from the same neuron over
%  different synapses.}

\section{The models for spiny dendrites}
\label{sec:models}
\subsection{Spike Diffuse Spike model}
\label{sec:SDS}

The SDS model \cite{Coombes1}, \cite{Coombes2} and \cite{Yulia}
describes a length of spiny dendritic by coupling a passive dendrite
to active spines. The cable is modelled by the passive cable equation,
and the spine head dynamics by the leaky integrate and fire model with
an imposed refractory time $\tau$.  The spines are attached at
discrete points along the cable by a spine stem with resistance, $r$.
The points at which the spines are attached can be spaced with any
spatial distribution but here we have chosen equally spaced
points. 

The membrane potential in the cable, V(x, t), is given by: 
\begin{eqnarray}
\label{eqn:SDScable}
\frac{\partial V}{\partial t} = D\frac{\partial^2 V}{\partial
  x^2}-\frac{V}{\tau} + {D r_a \rho(x)\frac{\hat{V}-V}{r}} + (\mu +
\nu g(V))\ast\frac{\partial Z(t)}{\partial t}
%\frac{\partial V}{\partial t} &=& D\frac{\partial^2 V}{\partial
%x^2}-\frac{V}{\tau} + {D r_a \rho(x)\frac{\hat{V}-V}{r}} + (\mu + \nu
%g(V))\ast\frac{\partial W(x, t)}{\partial t}\textrm{ ,} 
\end{eqnarray} 
with sealed end boundary conditions, $D = \frac{\lambda^2}{\tau}$ is
the diffusion coefficient, $\tau = r_mc_m$ is the membrane time
constant, $\lambda = \sqrt{\frac{d r_m}{4 r_a}}$ is the electronic
space constant, $r_a$ is the intracellular resistance per unit length,
$x\in[0, L]$ and $t\in[0, T]$, $L$ is the length of the cable and $T$
is the final time. We take $Z$ to either be a space-time Wiener
process $Z=W(x, t)$ chosen to satisfy the form of spatially correlated
noise we require or we take $Z(t)=K(t)$ a white noise or temporally
correlated noise path that is constant in space.  
We use $\ast$ to denote that the stochastic
integral can be interpreted in either an \ito or \strat sense and it
will be expressed by the usual convention: $\int g(V)dW$ or $\int
g(V)\circ dW$ respectively where appropriate.  
%and as such $\ast=\left \begin{array}{ll}
%.&\textrm{\ito calculus}\\
%\circ&\textrm{\strat calculus}
%\end{array}\right.$. 
$\rho(x) = \sum_{n\in \Gamma}\delta (x-x_n)$ is the density of spines,
attached at discrete points $x_n$, $\hat{V}(x_n, t)$ is the action
potential produced by the spine at the point $x_n$, the form of the
action potential is free to be chosen, and $r$ is the spine stem
resistance.  

The spine head dynamics are modelled by the stochastic
leaky integrate and fire (IF) model, upto firing threshold: 
\begin{eqnarray}
\label{eqn:SIF}
\frac{dU_n(t)}{dt} &=& \frac{V_n(t)}{\hat{C}r} - \epsilon U_n(t)
-\underbrace{\hat{C}h \sum_m \delta (t-T_n^m)}_{\textrm{Reset}}+ (\mu
+ \nu g(U_n(t)))\ast\frac{dZ_n(t)}{dt} .
%\frac{dU_n(t)}{dt} &=& \frac{V_n(t)}{\hat{C}r} - \epsilon U_n(t)
%-\underbrace{\hat{C}h \sum_m \delta (t-T_n^m)}_{\textrm{Reset}}+ (\mu
%+ \nu g(U_n(t)))\ast\frac{dW_n(x_n, t)}{dt}\textrm{ .}  
\end{eqnarray}
We take $Z$ to either be constant in space $Z_n(t)=K_n(t)$, where
$K_n$ is either a
white noise or temporally correlated path or $Z_n=W_n$ a spatially
correlated Wiener process evaluated at each spine at $x_n$. 
We choose the form of the 
multiplicative noise term to preserve the range of voltage $[0, 1]$:
$g(U) = U(1-U)$ while $U\in[0, 1]$ and zero otherwise,
\cite{Doering}. $\hat{C}$ is the capacitance, $r$ the spine stem
resistance and $\epsilon=\frac{1}{\hat{C}}\left(\frac{1}{r} +
  \frac{1}{\hat{r}}\right)$, with $\hat{r}$ as the membrane resistance
of the spine head. 
The m-th firing time of the nth spine, $T_n^m$, is governed by the
integrate and fire process: 
\begin{equation}
T_n^m = \inf\{t|U_n(t)\ge h, t>T_n^{m-1}+\tau_R\} \textrm{ ,}
\end{equation}
where $\tau_R$ is the refractory time period, during which the spine
is unable to fire.  This refractory time is introduced to mimic the
dynamics of the Hodgkin-Huxley (HH) model, which has a natural
refractory time.  At $U_n = h$ an action potential is injected into
the cable; the form of this injected potential can be chosen to be a
suitable function, in the SDS model described here, the function was
chosen to be a rectangular pulse, given by: 
\begin{equation}
\eta(t) = \eta_0\Theta (t)\Theta(\tau_s-t)\textrm{ ,}
\end{equation}
where $\Theta(t)$ is the Heaviside function, $\tau_s$ is the length of
time the pulse lasts for and $\eta_0$ is the strength (magnitude) of
the pulse.

These equations describing the dynamics of spiny dendritic tissue can
be solved using a combination of analytical and numerical techniques,
see \cite{Yulia} for a full description of the methodology.

The solution of this problem shows that the deterministic SDS model
supports the propagation of saltatory travelling waves along the
length of the cable \cite{Coombes2} and \cite{Yulia} provides a
preliminary investigation to these solutions with some noise in the
system but does not compare the speed of propagation or the effect of
spatial correlations. 

\subsection{Baer-Rinzel model}
\label{sec:BR}

The Baer-Rinzel (BR) model, \cite{Baer}, describes the voltage
evolution of a spiny dendrite using the Hodgkin Huxley (HH) equations
to describe the active properties of the spine heads: the voltage $U$
in the spine heads and the gating variables $m,n,h$.
They are coupled, with a certain density $\rho(x)$, to a uniform
passive cable, whose voltage, $V$, is modelled by the passive cable equation. 

We choose to include noise only to the m-dynamics since the other
gating variables are very sensitive to the noise and the m-gate is the
dominant variable as seen in reductions of the Hodgkin Huxley dynamics
to a two variable system, e.g. Fitzhugh Nagumo model, \cite{scott} and
\cite{Gerstner}. 
The equations for the stochastic BR model we consider are given by
three deterministic equations: 
\begin{eqnarray}
\hat{C}\frac{\partial U}{\partial t}&=&g_Kn^4(V_K-U) + g_{Na}hm^3(V_{Na}-U) + g_L(v_L-U) - \frac{U-V}{r}\\
\frac{dX}{dt} &=& \alpha_X(1-X) - \beta_XX\textrm{ ,}
\end{eqnarray}
where $X\in[n, h]$. Along with 
\begin{eqnarray}
\label{eqn:SBR}
C\frac{\partial V}{\partial t}&=&-g_L(V-V_L) + \frac{1}{r_a\pi
  d}\frac{\partial^2 V}{\partial x^2} + \rho(x)\frac{U-V}{r} + (\mu_c
+ \nu_cg_c(V))\ast\frac{\partial Z(x,t)}{\partial t}\\
\frac{dm}{dt} &=& \alpha_m(1-m) - \beta_mm + (\mu_m +
\nu_mg_m(m))\ast\frac{dZ_m}{dt}\textrm{ ,} 
\end{eqnarray}
%or for spatially correlated noise:
%\begin{eqnarray}
%\label{eqn:SBR2}
%C\frac{\partial V}{\partial t}&=&-g_L(V-V_L) + \frac{1}{r_a\pi
%d}\frac{\partial^2 V}{\partial x^2} + \rho(x)\frac{U-V}{r} + (\mu_c +
%\nu_cg_c(V))\ast\frac{\partial W(x, t)}{\partial t}\\ 
%\frac{dm}{dt} &=& \alpha_m(1-m) - \beta_mm + (\mu_m +
%\nu_mg_m(m))\ast\frac{dW_m}{dt}\textrm{ .} 
%\end{eqnarray}
the $\mu$'s give the strength of additive noise
and the $\nu$'s the strength of the multiplicative noise. 
We take either $Z(t)$=$K(t)$ is a constant in space and either white or
correlated in time or $Z(x,t)=W(x, t)$ for a spatially correlated Wiener
processes. 
Note that the equations for $X\in[m, n, h]$ and $U$ are coupled
together, at each point in space, by the cable or by the noise if it
is spatially correlated in the m-dynamics.  

When the noise is multiplicative (extrinsic), the stochastic integral
can be interpreted in either an \ito or \strat sense and 
we use the standard notation conventions: $\int g(V)dW$ for \ito and 
$\int g(V)\circ dW$ for \strat.
 The Stratonovich interpretation is often used in cases where the noise is fluctuating on a much faster scale than the system dynamics. The integrand $g(V)$ is evaluated at the mid-point of the interval i.e. $\frac{g(V(t)) + g(V(t+\Delta t))}{2}$, where $\Delta t$ is the time step and so the \strat interpretation can be thought of as averaging this fast behaviour in some way over the time interval.  However if the time-scales are much closer and the system responds to the noise on a similar time-scale then the \ito interpretation is more appropriate since it is non-anticipating, evaluating at the left hand end point of the time interval, \cite{Oskendal}, \cite{Kloeden}. \ito calculus is favoured in biological and financial applications due to this non-anticipating property i.e. only historical information is usually available. Despite this general assumptions each case should be considered individually, e.g. the results presented in \cite{Lindner03} show that simulations using a \strat interpretation are in closer agreement with the analytically predicted firing rates of a simple model of a type I neuron. Here we present results, \cite{thesis}, which differ depending on the noise interpretation used. It could be possible to use this discrepancy to determine the appropriate stochastic calculus for wave propagation in spiny dendrites.  
%{\bf Discuss choice. One interpretation of these two different
 % calculii is a time scale over which the noise arises. The \strat
  %interpretation often arises as the limit of a smooth apprximatoin
  %over each interval.{CHECK} An \ito interpretation is
  %non-anticipative and could be considered as the noise acting over a
  %small timesace. PLEASE CHECK ADD REFS ! We show below REF different
  %dynamics associated with the two interpretations - this suggests a
  %means to test wich interpretation OR WHY DO WE PICK ONE OR THE
  %OTHER.}
We choose the multiplicative noise functions $g$, in each of the
equation to ensure the fluctuations are added correctly to the resting
state of $V$ and $m$. With this in mind we choose the functions of the
multiplicative noise to be: $g_c(V)=-(65+V)$ and $g_m(m)=m(1-m)$ when
$m\in[0, 1]$ and zero otherwise. 

The spine density, $\rho(x)$ can be a constant as in the original BR
model \cite{Baer} which has been shown to support travelling waves, solitary,
multi-bump and periodic waves, \cite{Lord}. 

We extend the model here and consider here spatially 
dependent density, which has a parameter, $\kappa$, to control the
area of the spine stem attached to the cable. 
For spines centered at spatial points $x_n$: 
\begin{equation}
\label{eqn:rho_mol}
\rho(x) =  \sum_n \rho_{max} \xi_n \exp(-\kappa (x-x_n)^2)\textrm{ .}
\end{equation}
Here we have $\rho_{max}$ is the maximum value of the density, taken
to be the value used for the original BR model, and  
\begin{displaymath}
\xi_n = \left\{ \begin{array}{ll}
1 & \textrm{    if $x_n - \frac{d}{2}<x\le x_n + \frac{d}{2}$}\\
0 & \textrm{    otherwise,}
\end{array} \right.
\end{displaymath}
where $d$ is the spine spacing and $\kappa\in \mathbb{R}^+$ controls
the width of the spine stem. As $\kappa \to 0$, $\rho(x)\to
\rho_{max}$ and as $\kappa \to \infty$, $\rho(x) \to \sum_n \rho_{max}
\delta(x-x_n)$. Therefore at the $\kappa\to 0$ limit the model is the
BR model and at $\kappa \to\infty$ limit the model resembles the SDS
model with HH dynamics (instead of IF) in the spine heads. 

\subsection{Noise generation}
\label{sec:noise}
We use space/time white noise, temporally correlated and spatially correlated noise throughout the simulations to investigate the effects of noise on signal propagation in the previously described dendrite models. We consider a temporally correlated noise generated by the Ornstein-Uhlenbeck process, given as:
\begin{equation}
\label{eqn:OU2}
dK(t) = \beta(\theta - K(t))dt + \sigma db(t)\textrm{,}
\end{equation}
where $K(t)$ is a stochastic process called the Ornstein-Uhlenbeck
process, $\beta$ is a parameter which can adjust the time scale of the
correlation, called the mean reversion rate, $\theta$ is the mean to
which the process will revert to if given enough time, $\sigma$ is
another parameter which is called the volatility and $b$ is a Brownian
motion. To get a white noise path we simply set $\beta=1$ and
$\sigma=1$. 

The generation of a spatial correlated noise is introduced by a
process described in e.g. \cite{Tony} and \cite{Garcia}.  This form of
the spatially correlated noise is chosen to satisfy chosen properties
of the correlation length; we require that the correlation is over a
short range since we would not expect interference between distant
spines but would expect neighbouring spines to affect each other.   
We define a Q-Wiener process by the sum:
\begin{equation}
\label{eqn:corrsum}
W(x, t) = \sum_{j\ge 0}\lambda_j^{\frac{1}{2}}(x)\mathbf{e_j}(x)b_j(t)\textrm{ .}
\end{equation}
Assuming that $Q$ has the same eigenfunctions as the Laplacian (and
assuming Neumann boundary conditions on $[0, L]$), $\Delta$,
$\mathbf{e_j}(x) = \sqrt{\frac{2}{L}}\cos(\frac{\pi j x}{L})$ where
$j=1, 2, 3,\dots$ and $\mathbf{e_0}(x) = \sqrt{\frac{1}{L}}$,
$b_j(t)$ are standard Brownian motions and $\lambda_j$ are eigenvalues
of $Q$ which we choose here to satisfy a form of spatial correlation. 
We choose a correlation such that the noise is white in time and has
an exponential correlation in space of length is $\zeta$ (to satisfy
the assertion that nearby spines have more influence than distant).
The covariance and correlation function are given by and 
\begin{eqnarray}
\label{eqn:covar}
\mathbb{E}(W(x, t)W(x', t'))&=&F_c(x-x')\min\{t, t'\}\\
F_c(x) &=& \frac{1}{2\zeta}\exp(-\frac{\pi x^2}{4\zeta^2})\textrm{ .}
\label{eqn:corf}
\end{eqnarray}
Using this form of the correlation function we obtain \cite{thesis} the
following form for the eigenvalues of $Q$: $\lambda_j = \exp(-\frac{\pi j^2
  \zeta^2}{2L^2})$. 
If we require a noise which is white in space and time we can use
$Q=I$, this is non-trace class but has the eigenvalues $\lambda_j=1$,
$\forall j$ and as such the sum reduces to: $W(x,
t)=\sum_{j=1}^{\infty}e_j(x)b_j(t)$. The correlation length $\zeta$
can be chosen such that the strongest effect of the correlation is
felt over neighbouring spines and not, for example, the entire length
of the cable. Since we have chosen the form of the correlation to
ensure nearby spines have more influence on neighbours than distant
spines it would be counterproductive to use a long correlation
scale. We investigate by the correlation length the effect of
correlated signals into the dendrite as indicated in the introduction.

In order to simulate the different interpretations of the noise we require different numerical schemes. In the case of an \ito interpretation we use the Euler-Maruyama method and for the \strat interpretation we use the stochastic Heun method, for details on the implementation of these methods see \cite{thesis}.   

%{\bf Relate correlation length and spine spacing}

%{\bf Numerics : EM for \ito and Heun for \strat.}

\subsection{Speed of propagation of stochastic wave}
\label{sec:speed}

We rescale computed speeds by the deterministic speed,
$c_{det}=\frac{d}{t}$, where $d$ is the distance travelled along the
cable and $t$ is the time the wave takes to travel this distance; thus
the plotted speeds are given by $c=\frac{c_{noisy}}{c_{det}}$. This
rescaling allows comparison of the different models which have
different absolute values for the speed; it also makes it easy to see
in the graphs if the wave is speeding up ($c>1$) or slowing down
($c<1$), with respect to the deterministic wave speed, \cite{thesis}. 
In order to find the speed of any stochastic travelling wave
propagating on the cable we find the times, $t_1$ and $t_2$, at which
the wave crosses two points, $x_1$ and $x_2$, along the cable and use: 
\begin{equation}
  c_{noisy} = \frac{x_2 - x_1}{t_2 - t_1}\textrm{.}
\end{equation}
If the wave fails to reach $x_2$ then we say that the wave has failed
to propagate, \cite{thesis}. Failure is determined by the size of the
voltage in the cable at point $x_2$; if $V(x_2, t)\ge \theta$,
$\theta$ is a threshold value chosen from the voltage values in the
deterministic case, then it is still propagating. The threshold value
must be large enough that the voltage will only reach this level if
the voltage is close to that of the deterministic system and so we
avoid the case where the propagation is purely noise induced, i.e. we
are not measuring small fluctuations induced by noise only but we do
see the underlying signal too. We also impose the condition that the
wave must travel sequentially, i.e. each spine must fire in spatial
order from distal to proximal end of the dendrite. 

\section{Effects of the noise on speed of propagation}
\label{sec:effects}
Here we look at the results of the measuring the average speed of the
waves in the noisy BR and SDS models, over 100 realisations and for
the parameter values in \figref{fig:par}. We consider white noise,
temporally correlated noise and spatially correlated noise in both the
spine heads and the cable of each model (the noise is in the
m-dynamics for the BR model). We show error bars that represent the
standard deviation around the mean value of speed for each value of
noise intensity ($\pm
\textrm{S.D.(c)}=\sqrt{\mathbb{E}(c^2)-\left(\mathbb{E}(c)\right)^2}$). 
Table \ref{fig:par} contains the parameter values taken and their units.

%\subsection{Effect of multiplicative \ito noise in the BR and SDS models}
\subsection{Effect of multiplicative noise in the BR and SDS models}
\label{sec:Mresults}
The following results show the effect of multiplicative noise in both
models interpreted in the \ito sense. 
\begin{figure}[!htbp]
\begin{center}
\begin{tabular}{cc}
(a)&(b)\\
\resizebox{60mm}{!}{\includegraphics{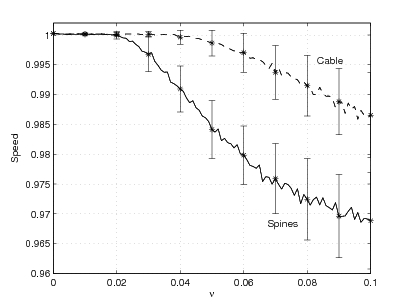}}&
\resizebox{60mm}{!}{\includegraphics{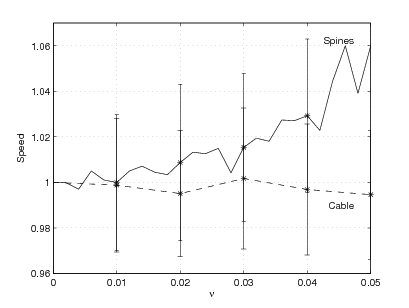}}\\
(c)&(d)\\
\resizebox{60mm}{!}{\includegraphics{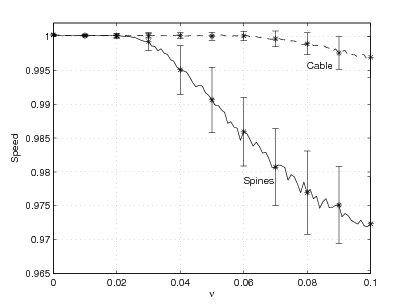}}&
\resizebox{60mm}{!}{\includegraphics{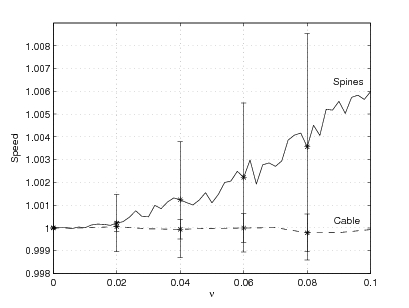}}
\end{tabular}
\caption{This figure shows how the speed of propagation changes as the
  noise intensity increases for both the SDS and BR models. Each plot
  shows the effect of noise in both the spines (solid lines) and the
  cable (dotted lines); the error bars are $\pm$ the standard
  deviation and get larger as the noise intensity increases. Plot (a)
  shows that space-time white noise in the spine and cable dynamics of the SDS
  model reduces the speed of propagation, much more so when the noise
  is in the spine heads. Plot (b) shows that white noise in the BR
  model increases the speed of propagation when the noise is driving
  the m-dynamics and that the cable is robust to this noise as the
  speed changes little with noise intensity. Plot (c) shows the SDS
  model driven by temporally correlated noise in the spine and cable dynamics and again
  the cable is more robust than the spines to the noise. As the noise
  intensity increases the speed of propagation decreases.  Plot (d)
  shows the speed of propagation in the BR model increasing with temporally correlated
  noise intensity when the noise is in the m-dynamics, again the cable
  seems robust to the effects of noise.}
\label{fig:speed}
\end{center}
\end{figure}
\Figref{fig:speed} shows that there is a difference in the behaviour
of the two models under the influence of synaptic noise
(i.e. multiplicative noise driving the spine dynamics); the speed of
propagation in the SDS model decreases as noise intensity increases
and the speed increases in the BR model with an increase in noise
intensity. In both cases the cable is more robust to the noise and
there is little or no variation in the speed as the noise intensity in
the cable increases. We can interpret this to mean that the synaptic
noise/noise in the input signal to be more important than background
noise in the propagation of signals in spiny dendrites. Since the
effect in the cable is negligible compared to the effect in the spines
we consider only the spines from now on. 

%%% \subsubsection{\ito vs Stratonovich}

\begin{figure}[!htbp]
\begin{center}
\begin{tabular}{cc}
(a)&(b)\\
\resizebox{60mm}{!}{\includegraphics{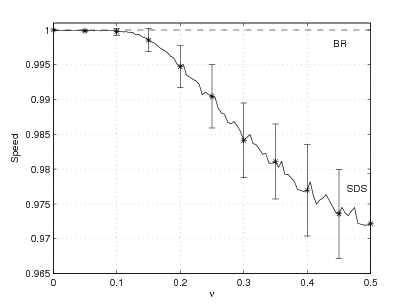}}&
\resizebox{60mm}{!}{\includegraphics{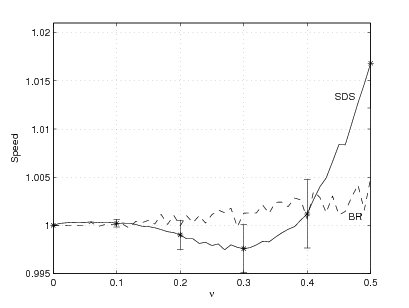}}\\
\end{tabular}
\caption{This figure shows the effect of spatially correlated noise in
  the spine dynamics of both models. Plot (a) shows effect of the \ito
  interpretation of the noise in the SDS model (solid line) and the BR
  model (dashed line). We observe a reduction in the speed as the
  noise intensity increases for the SDS model but no change for the BR
  model. Plot (b) shows the effect of the \strat interpretation of the
  noise in the SDS model (solid line) and BR model (dashed line). Both
  models show an increase in speed as the noise intensity increases,
  although the SDS model shows a larger increase.} 
\label{fig:BRx}
\end{center}
\end{figure}

We look at the effect of spatially correlated noise in the models,
again we look at the noise added to the m-dynamics for the BR model
and spine dynamics for the SDS model. \Figref{fig:BRx} shows the
effect of spatially correlated multiplicative noise on the speed of
propagation in the models as the noise intensity increases. We have
used both an \ito and Stratonovich interpretation to compare the
different noise models. The type of noise integral used in the
evaluation of the stochastic models has no qualitative effect on the
behaviour of the stochastic waves when the noise is white or OU but
when the noise is spatially correlated we observe a difference in the
behaviour for both models when the noise interpretation
changes. \figref{fig:BRx}, plot (a) shows, for \ito interpretation,
that the speed of a stochastic wave does not change, but for the SDS
model the speed decreases (as for white and OU noise). Plot (b) shows
an increase in speed for both the BR and SDS model, therefore the
choice of calculus when the noise is spatially correlated is
important, although it does not explain the difference in behaviour
for white or OU noise. 

\subsection{Small spatially correlated noise analysis of the BR model}
\label{sec:small_noise}

When we consider the noise to be interpreted by the \strat calculus we can implement a small noise expansion that exploits the non-zero mean of the noise term. 
This expansion of the stochastic BR model,
equations (5-7), to approximate the behaviour of the system without
having to simulate or solve any \strat stochastic integrals. We use
this approach for the original BR model since the constant density
makes the model amenable to analysis with the continuation package
AUTO-07P, \cite{auto}.  
This approach follows the working of \cite{Garcia} and \cite{Gardiner}
and the expansion is a standard perturbation approach. 
Applying this method to the stochastic BR model, \eqnref{eqn:SBR} and
converting to the travelling wave frame, using 
the standard anzatz $\xi = ct-x$, where $c$ is the wave speed, we
obtain: 
\begin{eqnarray}
\label{eqn:TW1}
V_\xi &=& \hat{W}\\
CcV_{\xi}&=&-g_L(V-V_L) + \frac{1}{r_a\pi d}V_{\xi} + \rho\frac{U-V}{r} - \nu_c^2F_c(0)(65+V)\\
c\hat{C}U_{\xi}&=&g_Kn^4(V_K-U)+g_{Na}hm^3(V_{Na}-U)+g_L(V_L-U)-\frac{U-V}{r}\\
cX_{\xi} &=& \alpha_X(1-X) - \beta_XX + \nu_X^2F_c(0)X(1-3X+2X^2)\textrm{,}
\label{eqn:TW2}
\end{eqnarray}
where $X\in[m, n, h]$.
We now have six coupled ODE's, instead of a PDE and four ODE's, and we
have new parameters, wave speed $c$ and noise intensity $\nu$. We now
want to investigate how these parameters effect the travelling wave
and also how the speed changes with noise intensity, to allow a
comparison with the full system. 
The continuation package AUTO-07P was used to investigate the 'new' system of equations in the travelling wave frame, similar to the work in \cite{Lord} for the deterministic BR model, and find how noise intensity changes the regions in parameter space where the solutions exist. 

\begin{figure}[!htbp]
\begin{center}
\begin{tabular}{cc}
(a)&(b)\\
\resizebox{50mm}{!}{\includegraphics{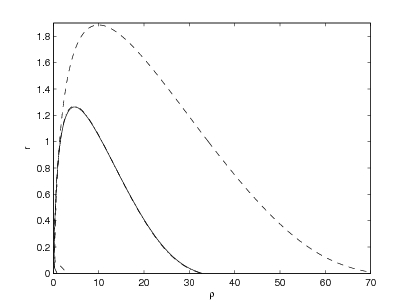}}&
\resizebox{50mm}{!}{\includegraphics{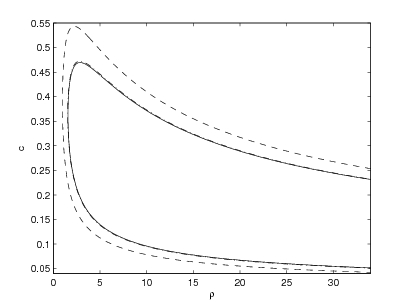}}\\
(c)&(d)\\
\resizebox{50mm}{!}{\includegraphics{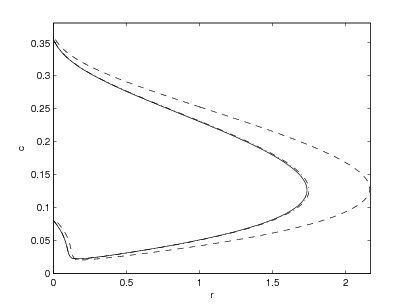}}
&\resizebox{50mm}{!}{\includegraphics{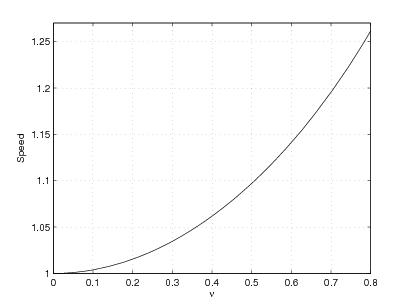}}
\end{tabular}
\caption{Examples of limit point diagrams for the deterministic and small noise case. In all plots the solid line is $\nu=0$, dot dashed line $\nu = 0.1$ and dashed line $\nu = 0.5$. Plot (a) is the spine density against the spine stem resistance, plot (b) shows spine density against the speed of the wave and plot (c) shows the spine stem resistance against the speed of the wave. Each plot show an increase in the area of existence in parameter space for the travelling waves as the noise intensity increases. Plot (d) shows the speed of the wave as the noise intensity in the spine head dynamics is increased.}
\label{fig:auto}
\end{center}
\end{figure}

We seek localised solutions of \eqnref{eqn:TW1} to \eqnref{eqn:TW2} as a homoclinic orbit which is approximated by a large periodic orbit (T=100). We can then continue in various combinations of parameters to find the regions of the parameter space where the travelling waves exist. We are interested in $\rho$ the spine density, $r$ the spine stem resistance, and $c$ the wave speed. To find limits for existence we investigate the effect of the noise intensity. In \figref{fig:auto} we show the areas of existence of the travelling wave in the $\rho$-$r$ parameter space, plot (a), and the limit point diagram for $\rho$ and the speed of the waves $c$ in plot (b); finally plot (c) shows the existence of travelling waves in $r$, $c$ parameter space. The fast branch of this figure (plot (b)) is the stable branch and shows the waves we observe in the direct simulations. 

In \figref{fig:auto} we look at $\nu_m=0$ (solid line), the deterministic case, $\nu_m=0.1$ (dot-dashed line) and $\nu_m=0.5$ (dashed line). In all three plots the area in which the waves exist in parameter space increases as the noise intensity increases. We see from this figure, plot (d) that the continuation results agree with the results for direct simulation of the BR model in the \strat sense, \figref{fig:BRx}.

To investigate the effect of the noise on the existence of these waves in parameter space, we choose a point on the fast branch of the deterministic $\rho$-c curve and continue from there in two parameters, for example $\nu_m$ (spatially correlated noise intensity in the m-dynamics) and $c$. We can then see how the wave speed changes as noise is added into the equations, \figref{fig:auto}, plot (d).  

\subsection{Baer - Rinzel model with variable density, $\rho(x)$}
\label{sec:BR_SDS}
We now look at the speed of a noisy wave as the parameter $\kappa$
changes, and so as the spine stem width changes, for an explanation of
the contradictory results from the BR and SDS models with white and OU
noise. From the behaviour we have observed so far we expect as
$\kappa$ increases and so the model changes from a BR type, $\kappa =
0$, to an SDS type, $\kappa = 1000$, then the plot of the noisy wave
speed should cross the plot of the deterministic wave speed, since the
BR model speeds up with the inclusion of noise and we expect the BR
model with discrete spine density (as in SDS model) to decrease in
speed if the dynamics are not as important as the spine stem
morphology. 
First we consider the new BR model without noise, as the spine stem
changes through the increase of parameter $\kappa$. 
\begin{figure}[h]
\begin{center}
\begin{tabular}{cc}
(a)&(b)\\
\resizebox{60mm}{!}{\includegraphics{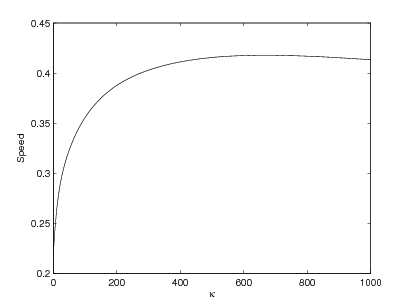}}&
\resizebox{60mm}{!}{\includegraphics{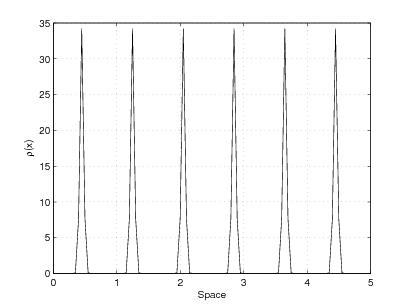}}
\end{tabular}
\caption{Plot (a) shows the deterministic speed of waves in the
  modified BR model as the parameter $\kappa$ increases from the BR
  limit ($\kappa=0$) to the SDS limit ($\kappa \to \infty$). There is
  an optimal value of $\kappa$ which maximises the speed of the wave,
  this occurs at $\kappa=670$. Plot (b) shows an examples of the space dependent density. Here $\kappa=670$ and $d=0.8$ is the spine spacing (distance between peaks).} 
\label{fig:Ninc_determ}
\end{center}
\end{figure}
\Figref{fig:Ninc_determ}, plot (a), shows the speed of a deterministic wave in the modified BR model as $\kappa$ changes and so the spine stem changes from a continuum like the original BR model to a discrete distribution of spines as in the SDS model. There is an optimal value of $\kappa$ which maximises the deterministic wave speed. This shows that the distribution or spine stem morphology is of importance in the propagation of action potentials, it can alter the speed of the AP on the dendrite. Plot (b) shows an example of the spine density for $\kappa=670$. The spatial discretisation will effect the density since a larger $\Delta x$ will approximate a discrete stem at a smaller value of $\kappa$ than a smaller step. The value of speed measured as $\Delta x$ changes is small and the overall behaviour (as $\kappa$ changes) remains the same.

\begin{figure}[h]
\begin{center}
\begin{tabular}{cc}
(a)&(b)\\
\resizebox{50mm}{!}{\includegraphics{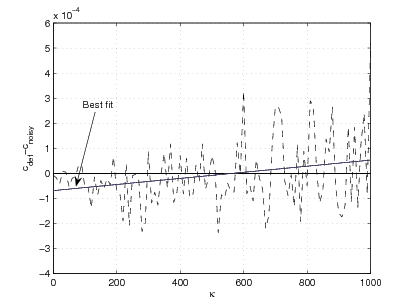}}&
\resizebox{50mm}{!}{\includegraphics{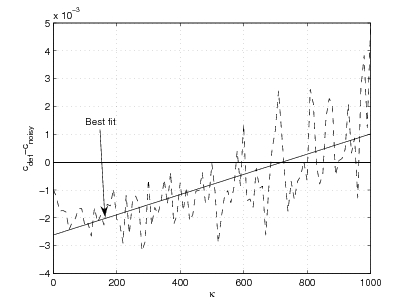}}
\end{tabular}
\caption{This figure shows the difference between the deterministic
  speed and the speed of the noisy wave as $\kappa$ increases. The noise is interpreted in the \ito sense. Plot
  (a) shows $c_{det}-c_{noisy}$, as $\kappa$ increases, where the
  noise intensity is $\nu=0.02$ and plot (b) has noise intensity
  $\nu=0.15$. The dotted lines show the difference in the speed,
  starting negative as the noisy speed is faster than the
  deterministic wave speed and becoming positive as the noisy wave
  becomes slower than the deterministic wave. In plot (a) the
  difference in wave speed is very small but as the noise intensity
  increases, plot (b), the difference increases and the dotted line is
  above the line $y=0$ for a larger range of $\kappa$.
} 
\label{fig:Ninc}
\end{center}
\end{figure}

We now consider this variable density BR model with different levels of noise in the spine head dynamics (m-dynamics), in the \ito sense, as $\kappa$ changes. We plot the difference $c_{det}-c_{noisy}$, to show how the noisy wave changes with respect to the deterministic wave as $\kappa$ increases. In order to make it clearer we also show the linear best fit line. It is clear in \figref{fig:Ninc} that the trend is for the wave to be slower ($c_{det}-c_{noisy}>0$) than the deterministic value at the SDS limit (large $\kappa$) and faster ($c_{det}-c_{noisy}<0$) at the BR limit ($\kappa=0$). \figref{fig:Ninc} plot (a) shows the noise intensity at $\nu=0.02$ and plot (b) at $\nu=0.15$.
It is clear that as the parameter $\kappa$ increases the behaviour of
the noisy wave in the new BR model also changes. When $\kappa=0$ the
model is the original BR model and the stochastic wave is, on average,
faster than the deterministic wave, and as $\kappa$ reaches large
values then the wave speed decreases below that of the deterministic
wave, mimicking the behaviour of a stochastic wave in the SDS
model. This goes to show that the form of the density has an important
effect on the behaviour of the model. 

\section{Additive noise}
\label{sec:add}
\begin{figure}[!htbp]
\begin{center}
\begin{tabular}{cc}
(a)&(b)\\
\resizebox{50mm}{!}{\includegraphics{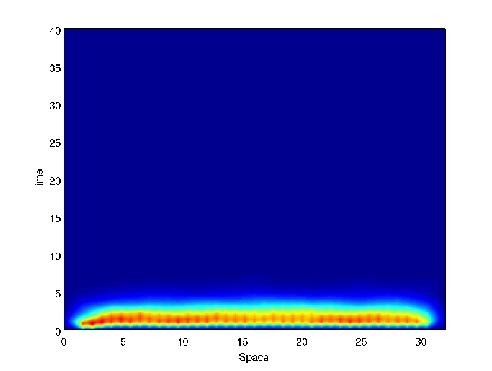}}&
\resizebox{50mm}{!}{\includegraphics{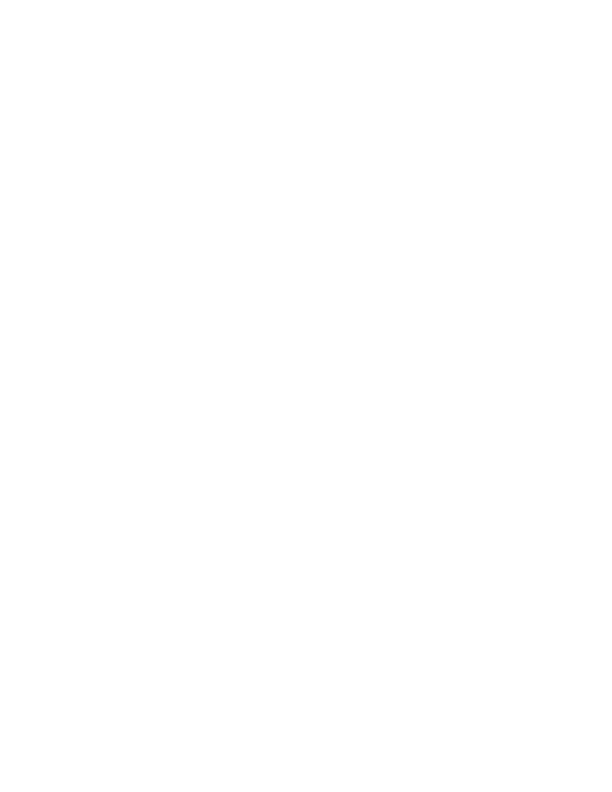}}
\end{tabular}
\end{center}
\caption{This figure shows the mean voltage in the cable for both the SDS (plot (a)) and BR (plot (b)) model. The colour shows the value of the voltage for each point in space (x-axis) and time (y-axis). As can be seen in both plots the wave is near horizontal which shows that all spines are firing at the same time. Plot (a) shows only one wave due to the long refractory time in the SDS model (which can be controlled by parameter $\tau$) and Plot (b) has many waves due to the natural refractory time of the BR model. }
\label{fig:Veg}
\end{figure}

So far we have only considered multiplicative noise and now briefly show that additive noise in the models can have some interesting effects. When the noise is additive white in the spine heads of both models we can observe synchrony for small levels of noise intensity, see \figref{fig:Veg}.

\begin{figure}[!htbp]
\begin{center}
\begin{tabular}{cc}
(a)&(b)\\
\resizebox{50mm}{!}{\includegraphics{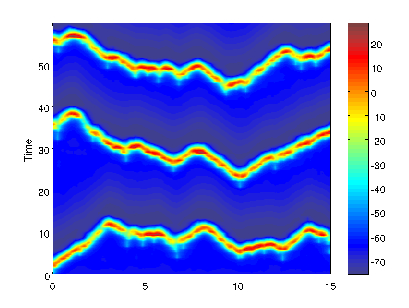}}&
\resizebox{50mm}{!}{\includegraphics{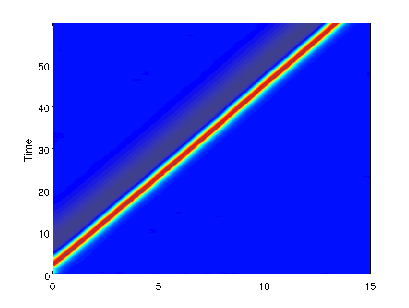}}
\end{tabular}
\caption{This figure shows a sample of the voltage in the cable when there is noise present in the spine head m-dynamics. Plot (a) has additive white noise of strength $\mu=0.01$, plot (b) has additive OU noise strength $\mu=0.01$ and $\beta=2$. It can be seen that the waves travel out of order when the noise is white but when a temporal correlation is added the waves regain their sequential travel.}
\label{fig:Vspines}
\end{center}
\end{figure}

The OU noise (additive and multiplicative) in the m-dynamics of the BR model helps to stabilise waves which were out of order in the white noise case. \Figref{fig:Vspines} plots (a) and (b) show the voltage in the cable with additive noise in the m-dynamics, $\mu=0.01$, for white and OU noise respectively. It is clear that the temporal correlation of $\beta=2$ promotes a travelling wave; it could do this by matching some internal time scale in the BR model. Plots (c) and (d) show the voltage in the cable when the noise is multiplicative in the m-dynamics, $\nu=0.1$, for white and $\beta=2$ correlated noise; again the correlation promotes a sequential travelling wave. 

\begin{figure}[!htbp]
\begin{center}
\begin{tabular}{cc}
(a)&(b)\\
\resizebox{50mm}{!}{\includegraphics{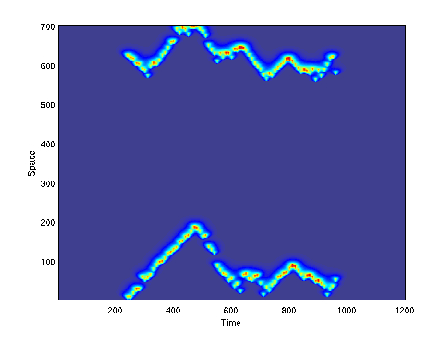}}&
\resizebox{50mm}{!}{\includegraphics{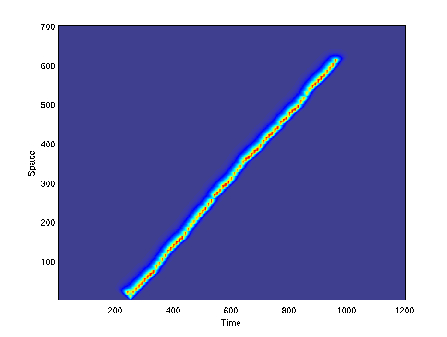}}
\end{tabular}
\end{center}
\caption{We show here the effect of the correlation length of additive noise in the spine heads of the SDS model for a fixed level of noise intensity. Plot (a) shows a short correlation scale and plot (b) a longer scale and it is clear that the longer correlation scale restores the sequential firing which is lost when the noise has a short correlation length. }
\label{fig:corrlen}
\end{figure}

The additive spatially correlated noise in the SDS model can stabalise non-sequential travelling waves as the correlation scale increases. \Figref{fig:corrlen} shows additive spatially correlated noise in the spines of the SDS model; plot (a) shows the voltage in the cable when the spines are coupled by a short correlation length and the wave fires out of order and plot (b) shows the voltage when the correlation length is longer and this restores the sequential firing of the spines and so the wave travels smoothly again.
The additive spatially correlated noise in the SDS model can stabalise non-sequential travelling waves as the correlation scale increases. \Figref{fig:corrlen} shows additive spatially correlated noise in the spines of the SDS model; plot (a) shows the voltage in the cable when the spines are coupled by a short correlation length and the wave fires out of order and plot (b) shows the voltage when the correlation length is longer and this restores the sequential firing of the spines and so the wave travels smoothly again.

\section{Discussion}
\label{sec:discuss}

We set out to investigate the effects of different types of noise on models of spiny dendrites and have provided a comparison of white, temporally and spatially correlated noise in the spike-diffuse-spike model and the Baer and Rinzel model. In all cases noise in the cable equation has little effect on the speed of propagation which suggests that the cable is robust whereas the spine dynamics are more sensitive to noisy input. This makes sense with regards to the structure of the models since the cable is passive and diffusive as opposed to the active properties of the spine heads.

When the noise is in the SDS model the speed of any travelling wave decreases as the noise intensity increases but when the noise is in the BR model the speed of travelling waves increases as the noise intensity increases, when the noise is in the spines. The main differences between the two models are the dynamics used to describe the evolution of the spine head voltage (integrate and fire dynamics in the SDS model and the Hodgkin Huxley equations in the BR model) and the spine density, $\rho(x)$ (discretely attached equally spaced spines in the SDS model and a constant in the BR model). This difference could be used to decide which model is a more accurate description of the real dendrite if an experiment could be devised in which the speed of an injected pulse travels the length of a dendrite with noise present. In an attempt to discover which of the differences in the models produced the difference in behaviour we investigated the BR model with a spatially dependent density. In the 'SDS' limit where the spines are attached discretely the model can be thought to be the SDS model with Hodgkin Huxley dynamics in the spine heads and this model does act like the SDS model when we add noise to the system. Therefore we can conclude that the type of spine dynamics is not the reason for the difference in behaviour between the SDS and BR models and we investigate the affect of the spine density and so the spine stem. As we change the spine density from one limit to the other the cases in between are like looking at spines attached to the cable with a spine stem that has an area of attachment to the cable. The speed of the wave in the deterministic version of this model  increases from the BR limit to the SDS limit and there is an optimal value for the spine attachment that gives a maximum value of the speed. When there is noise in the system and we are changing the parameter from BR to SDS limits, the average speed of the wave starts faster than the deterministic wave speed and decreases below the deterministic speed at the SDS limit, in agreement with the previous results from the original SDS and BR models. Although this spatially dependent spine density gives some insight to the importance of the spine stem the model could be improved by a more realistic physical description of the way the spine stems are attached to the cable and perhaps by introducing some random distribution of the spines to investigate how this affects the behaviour. To take this a step further there could even be a simple learning rule introduced which could move the spines in relation to the levels of activity and give a time dependent spine distribution; \cite{Verzi2} looks at a simple activity dependent spine plasticity in the BR model, so perhaps this could be extended to include noise.
When both models are influenced by spatially correlated noise the interpretation of the noise (\ito or \strat) also changes the behaviour of the travelling waves. The Baer and Rinzel model is not affected when the \ito interpretation is used but the \strat interpretation causes an increase in speed as seen in direct simulation and through a small noise analysis, \secref{sec:small_noise}. Likewise the SDS model shows a change in behaviour: speed reduction in the \ito case and an increase in the \strat case. This begs the question which is the correct noise interpretation for this situation? Is there some way to investigate, experimentally, which is the better noise interpretation for a length of dendrite?  
 
Additive noise in both the dendrite models can induce synchronous behaviour in the spines. For very small levels of noise the systems are fairly robust to the noise; all waves fully propagate and are only subject to a small change in the speed. As the strength of additive noise increases the spines in the models begin to fire out of order and seemingly in a random fashion until a synchronous behaviour takes over and they all fire simultaneously, similarly observed in \cite{Newhall} and \cite{Newhall2}. When the noise in the SDS model is additive and spatially correlated in the spine heads then the correlation scale can play a role in restoring sequential firing that has been destroyed by the noise, e.g. when the noise intensity is fixed a short correlation scale displays the out of order firing but as the correlation scale is increased the wave travels in a sequential fashion. 

These results suggest that there is an optimal spine density for speed of propagation in dendrites; too sparse and propagation slows down and ultimately stops, too many and we again see a reduction in speed. A spatially correlated noise can also enhance propagation, we used a correlation length of 3 spine spacings which suggests that a localised input to the dendrite, at synapses on neighbouring spines, can increase propagation speed. There are obvious places to extend this work; computationally further investigation of the optimal correlation scale, and propagation of a noisy signal on a dendrite with a random spine distribution. Experimentally it would be useful to verify which noise interpretation is more realistic and if noise speeds up or slows down a real travelling wave.

%\begin{figure}
  %\includegraphics{example.eps}
%\end{figure}
%\begin{figure*}
  %\includegraphics[width=0.75\textwidth]{example.eps}
%\end{figure*}
%\begin{table}
%\begin{tabular}{lll}
%\hline\noalign{\smallskip}
%\end{tabular}
%\end{table}

%\begin{acknowledgements}
%\end{acknowledgements}

\begin{appendix}
Table \ref{fig:par} gives the parameter values used throughout this work and
their respective units. The values are biologically realistic, and the bracketed values are the non-dimensional values
used.  The number of spines used in the SDS model was chosen such that they were equally spaced along the (non-dimensional) length of cable used, and such that wave propagation was present.

\begin{center}
\begin{tabular}{llll}
\label{fig:par}
\textbf{Symbol}&\textbf{Name}&\textbf{Value}&\textbf{Unit}\\
\hline\noalign{\smallskip}
$V$&Cable Voltage&-&$mV$\\
$U$&Spine head Voltage&-&$mV$\\
$R_m$&Transmembrane Resistance&2500 (1)&$\Omega cm^2$\\
$R_a$&Intracellular Resistance&70 (1)&$\Omega cm$\\
$C_m$&Transmembrane Capacitance&1 (1)&$\mu F cm^{-2}$\\
$\hat{C}$&Transmembrane Capacitance of Spine head&1 (1)&$\mu F cm^{-2}$\\
$\hat{r}$&Transmemrane Resistance of Spine head&2500 (1)&$\Omega cm^2$\\
$\mu$&Strength of additive noise&-&-\\
$\nu$&Strength of multiplicative noise&-&-\\
$a$&Dendritic Diameter&0.36 (1)&$\mu m$\\
$\lambda=\sqrt{aR_m}{4R_a}$&Electronic Length Scale& (1)&-\\
$\tau=R_mC_m$&Electonic time constant& (1)&-\\
$D=\frac{\lambda^2}{\tau}$&Diffusion Coefficient&(1)&-\\
$\tau_R$&Refractory time&-&-\\
$L$&Length of dendrite&1-2&$mm$\\
$N_{spines}$&Number of spines in SDS model&81&-\\
$\tau_s$&Length of time pulse lasts in SDS model&-&-\\
$h$&Voltage threshold in spine head for SDS model&0.04&-\\
$m$&Sodium activation particle&-&-\\
$h$&Sodium inactivation particle&-&-\\
$n$&Potassium activation particle&-&-\\
$G_{Na}$&Maximum sodium conductance&120&$mScm^{-2}$\\
$G_{K}$&Maximum potassium conductance&36&$mScm^{-2}$\\
$G_{L}$&Maximum leakage conductance&0.3&$mScm^{-2}$\\
$V_{Na}$&Sodium reversal potential&50&$mV$\\
$V_{K}$&Potassium reversal potential&-77&$mV$\\
$V_{L}$&Leakage reversal potential&-54.402&$mV$\\
$\rho$&Spine density&-&-\\
$d$&Spine spacing&(0.8 or 1)&$\mu m$\\
$\eta_0$&Strength of action potential (AP) pulse in SDS model&-&-\\
\end{tabular}
\end{center}

\end{appendix}


%\begin{thebibliography}{}
% etc
%\end{thebibliography}

\bibliographystyle{plain}
\bibliography{refs.bib}
\end{document}